\documentclass[11pt]{article}
\usepackage[caption = false]{subfig}

\usepackage{amsmath}
\usepackage{amsfonts}
\usepackage[mathscr]{eucal}
\usepackage{amsthm}
\usepackage{amssymb}
\usepackage[noadjust]{cite}
\usepackage{enumerate}
\usepackage{hyperref}
\usepackage{graphicx}

\usepackage{color}
\usepackage{soul}

\graphicspath{{./}{figs/}}

\DeclareMathOperator{\diam}{diam}

\newcommand{\cT}{\mathcal{T}}

\begin{document}

\newcommand{\rf}[1]{(\ref{#1})}
\newcommand{\mmbox}[1]{\fbox{\ensuremath{\displaystyle{ #1 }}}}	%

\newcommand{\hs}[1]{\hspace{#1mm}}
\newcommand{\vs}[1]{\vspace{#1mm}}

\newcommand{\ri}{{\mathrm{i}}}
\newcommand{\re}{{\mathrm{e}}}
\newcommand{\rd}{\mathrm{d}}

\newcommand{\R}{\mathbb{R}}
\newcommand{\Q}{\mathbb{Q}}
\newcommand{\N}{\mathbb{N}}
\newcommand{\Z}{\mathbb{Z}}
\newcommand{\C}{\mathbb{C}}
\newcommand{\K}{{\mathbb{K}}}

\newcommand{\cA}{\mathcal{A}}
\newcommand{\cB}{\mathcal{B}}
\newcommand{\cC}{\mathcal{C}}
\newcommand{\cS}{\mathcal{S}}
\newcommand{\cD}{\mathcal{D}}
\newcommand{\cH}{\mathcal{H}}
\newcommand{\cI}{\mathcal{I}}
\newcommand{\cItilde}{\tilde{\mathcal{I}}}
\newcommand{\cIhat}{\hat{\mathcal{I}}}
\newcommand{\cIcheck}{\check{\mathcal{I}}}
\newcommand{\cIstar}{{\mathcal{I}^*}}
\newcommand{\cJ}{\mathcal{J}}
\newcommand{\cM}{\mathcal{M}}
\newcommand{\cP}{\mathcal{P}}
\newcommand{\cV}{{\mathcal V}}
\newcommand{\cW}{{\mathcal W}}
\newcommand{\scrD}{\mathscr{D}}
\newcommand{\scrS}{\mathscr{S}}
\newcommand{\scrJ}{\mathscr{J}}
\newcommand{\sD}{\mathsf{D}}
\newcommand{\sN}{\mathsf{N}}
\newcommand{\sS}{\mathsf{S}}
 \newcommand{\sT}{\mathsf{T}}
 \newcommand{\sH}{\mathsf{H}}
 \newcommand{\sI}{\mathsf{I}}
 
\newcommand{\bs}[1]{\mathbf{#1}}
\newcommand{\bb}{\mathbf{b}}
\newcommand{\bd}{\mathbf{d}}
\newcommand{\bn}{\mathbf{n}}
\newcommand{\bp}{\mathbf{p}}
\newcommand{\bP}{\mathbf{P}}
\newcommand{\bv}{\mathbf{v}}
\newcommand{\bx}{\mathbf{x}}
\newcommand{\by}{\mathbf{y}}
\newcommand{\bz}{{\mathbf{z}}}
\newcommand{\bxi}{\boldsymbol{\xi}}
\newcommand{\boldeta}{\boldsymbol{\eta}}	%

\newcommand{\ts}{\tilde{s}}
\newcommand{\tGamma}{{\tilde{\Gamma}}}
 \newcommand{\tbx}{\tilde{\bx}}
 \newcommand{\tbd}{\tilde{\bd}}
 \newcommand{\txi}{\xi}
 
\newcommand{\done}[2]{\dfrac{d {#1}}{d {#2}}}
\newcommand{\donet}[2]{\frac{d {#1}}{d {#2}}}
\newcommand{\pdone}[2]{\dfrac{\partial {#1}}{\partial {#2}}}
\newcommand{\pdonet}[2]{\frac{\partial {#1}}{\partial {#2}}}
\newcommand{\pdonetext}[2]{\partial {#1}/\partial {#2}}
\newcommand{\pdtwo}[2]{\dfrac{\partial^2 {#1}}{\partial {#2}^2}}
\newcommand{\pdtwot}[2]{\frac{\partial^2 {#1}}{\partial {#2}^2}}
\newcommand{\pdtwomix}[3]{\dfrac{\partial^2 {#1}}{\partial {#2}\partial {#3}}}
\newcommand{\pdtwomixt}[3]{\frac{\partial^2 {#1}}{\partial {#2}\partial {#3}}}
\newcommand{\bnabla}{\boldsymbol{\nabla}}
\newcommand{\dive}{\boldsymbol{\nabla}\cdot}
\newcommand{\curl}{\boldsymbol{\nabla}\times}
\newcommand{\Phixy}{\Phi(\bx,\by)}
\newcommand{\PhiOxy}{\Phi_0(\bx,\by)}
\newcommand{\dxPhixy}{\pdone{\Phi}{n(\bx)}(\bx,\by)}
\newcommand{\dyPhixy}{\pdone{\Phi}{n(\by)}(\bx,\by)}
\newcommand{\dxPhiOxy}{\pdone{\Phi_0}{n(\bx)}(\bx,\by)}
\newcommand{\dyPhiOxy}{\pdone{\Phi_0}{n(\by)}(\bx,\by)}

\newcommand{\eps}{\varepsilon}
\newcommand{\real}[1]{{\rm Re}\left[#1\right]} %
\newcommand{\im}[1]{{\rm Im}\left[#1\right]}
\newcommand{\ol}[1]{\overline{#1}}
\newcommand{\ord}[1]{\mathcal{O}\left(#1\right)}
\newcommand{\oord}[1]{o\left(#1\right)}
\newcommand{\Ord}[1]{\Theta\left(#1\right)}

\newcommand{\hsnorm}[1]{||#1||_{H^{s}(\bs{R})}}
\newcommand{\hnorm}[1]{||#1||_{\tilde{H}^{-1/2}((0,1))}}
\newcommand{\norm}[2]{\left\|#1\right\|_{#2}}
\newcommand{\normt}[2]{\|#1\|_{#2}}
\newcommand{\on}[1]{\Vert{#1} \Vert_{1}}
\newcommand{\tn}[1]{\Vert{#1} \Vert_{2}}

\newcommand{\xt}{\mathbf{x},t}
\newcommand{\PhiF}{\Phi_{\rm freq}}
\newcommand{\cone}{{c_{j}^\pm}}
\newcommand{\ctwo}{{c_{2,j}^\pm}}
\newcommand{\cthree}{{c_{3,j}^\pm}}

\newtheorem{thm}{Theorem}[section]
\newtheorem{lem}[thm]{Lemma}
\newtheorem{defn}[thm]{Definition}
\newtheorem{prop}[thm]{Proposition}
\newtheorem{cor}[thm]{Corollary}
\newtheorem{rem}[thm]{Remark}
\newtheorem{conj}[thm]{Conjecture}
\newtheorem{ass}[thm]{Assumption}
\newtheorem{example}[thm]{Example} %

\newcommand{\tH}{\widetilde{H}}
\newcommand{\Hze}{H_{\rm ze}} 	%
\newcommand{\uze}{u_{\rm ze}}		%
\newcommand{\dimH}{{\rm dim_H}}
\newcommand{\dimB}{{\rm dim_B}}
\newcommand{\IntClosOm}{\mathrm{int}(\overline{\Omega})}
\newcommand{\IntClosOmOne}{\mathrm{int}(\overline{\Omega_1})}
\newcommand{\IntClosOmTwo}{\mathrm{int}(\overline{\Omega_2})}
\newcommand{\Ccomp}{C^{\rm comp}}
\newcommand{\tCcomp}{\tilde{C}^{\rm comp}}
\newcommand{\uC}{\underline{C}}
\newcommand{\utC}{\underline{\tilde{C}}}
\newcommand{\oC}{\overline{C}}
\newcommand{\otC}{\overline{\tilde{C}}}
\newcommand{\capcomp}{{\rm cap}^{\rm comp}}
\newcommand{\Capcomp}{{\rm Cap}^{\rm comp}}
\newcommand{\tcapcomp}{\widetilde{{\rm cap}}^{\rm comp}}
\newcommand{\tCapcomp}{\widetilde{{\rm Cap}}^{\rm comp}}
\newcommand{\hcapcomp}{\widehat{{\rm cap}}^{\rm comp}}
\newcommand{\hCapcomp}{\widehat{{\rm Cap}}^{\rm comp}}
\newcommand{\tcap}{\widetilde{{\rm cap}}}
\newcommand{\tCap}{\widetilde{{\rm Cap}}}
\newcommand{\ccap}{{\rm cap}}
\newcommand{\ucap}{\underline{\rm cap}}
\newcommand{\uCap}{\underline{\rm Cap}}
\newcommand{\cCap}{{\rm Cap}}
\newcommand{\ocap}{\overline{\rm cap}}
\newcommand{\oCap}{\overline{\rm Cap}}
\DeclareRobustCommand
{\mathringbig}[1]{\accentset{\smash{\raisebox{-0.1ex}{$\scriptstyle\circ$}}}{#1}\rule{0pt}{2.3ex}}
\newcommand{\cirH}{\mathringbig{H}}
\newcommand{\cirHs}{\mathringbig{H}{}^s}
\newcommand{\cirHt}{\mathringbig{H}{}^t}
\newcommand{\cirHm}{\mathringbig{H}{}^m}
\newcommand{\cirHzero}{\mathringbig{H}{}^0}
\newcommand{\deO}{{\partial\Omega}}
\newcommand{\OO}{{(\Omega)}}
\newcommand{\Rn}{{(\R^n)}}
\newcommand{\Id}{{\mathrm{Id}}}
\newcommand{\gap}{\mathrm{Gap}}
\newcommand{\ggap}{\mathrm{gap}}
\newcommand{\isom}{{\xrightarrow{\sim}}}
\newcommand{\half}{{1/2}}
\newcommand{\mhalf}{{-1/2}}
\newcommand{\inter}{{\mathrm{int}}}

\newcommand{\Hsp}{H^{s,p}}
\newcommand{\Htq}{H^{t,q}}
\newcommand{\tHsp}{{{\widetilde H}^{s,p}}}
\newcommand{\SP}{\ensuremath{(s,p)}}
\newcommand{\Xsp}{X^{s,p}}

\newcommand{\dd}{{d}}\newcommand{\pp}{{p_*}}

\newcommand{\Rnn}{\R^{n_1+n_2}}
\newcommand{\Tr}{{\mathrm{Tr}}}

\renewcommand{\arraystretch}{1.7}
\renewcommand{\tH}{\widetilde{H}{}}
\allowdisplaybreaks[4]
\newcommand{\cK}{\mathcal{K}}

\title{
Discontinuous piecewise polynomial approximation 
on 
non-Lipschitz
domains
}
\author{
David P. Hewett$^{\text{a}}$\\
$^{\text{a}}${\footnotesize Department of Mathematics, University College London, London, United Kingdom. \texttt{d.hewett@ucl.ac.uk}
}}

\maketitle
\renewcommand{\thefootnote}{\arabic{footnote}}

\begin{abstract}
We prove best approximation error estimates for discontinuous piecewise polynomial approximation in fractional Sobolev spaces on non-Lipschitz meshes of non-Lipschitz domains. In particular, the boundary of the domain, and the boundaries of the mesh elements, can be fractal. 
\end{abstract}


\section{Introduction}

Piecewise polynomial approximation on 
a domain 
$\Omega\subset\R^n$ %
is central to many numerical methods, including the finite element method (FEM) \cite{brenner2008mathematical,Ciarlet78} and integral equation methods \cite{Kress,atkinson1997numerical} such as the boundary element method \cite{sauter-schwab11}. 
In conforming FEM, the piecewise polynomial approximation space must be globally continuous. This imposes strict constraints on the mesh geometry, and classically one tends to work with conforming simplicial meshes. 
However, in discontinuous Galerkin (dG) FEM, and in many integral equation methods, one works with globally discontinuous piecewise polynomials, which introduces considerably more flexibility in the type of mesh that can be used. When the underlying domain 
has a complicated structure, the use of 
meshes with 
non-simplicial and even non-polytopal elements offers significantly more efficient and accurate approximation of the domain geometry than can be achieved using traditional simplicial meshes. 
In dG-FEM this idea is well established, with the theoretical state of the art being \cite{cangiani2022hp}, which proves optimal $hp$ best approximation error estimates 
for Lipschitz domains and a class of meshes with Lipschitz elements.

\begin{figure}
\centering
\subfloat[]{
\includegraphics[height=50mm]{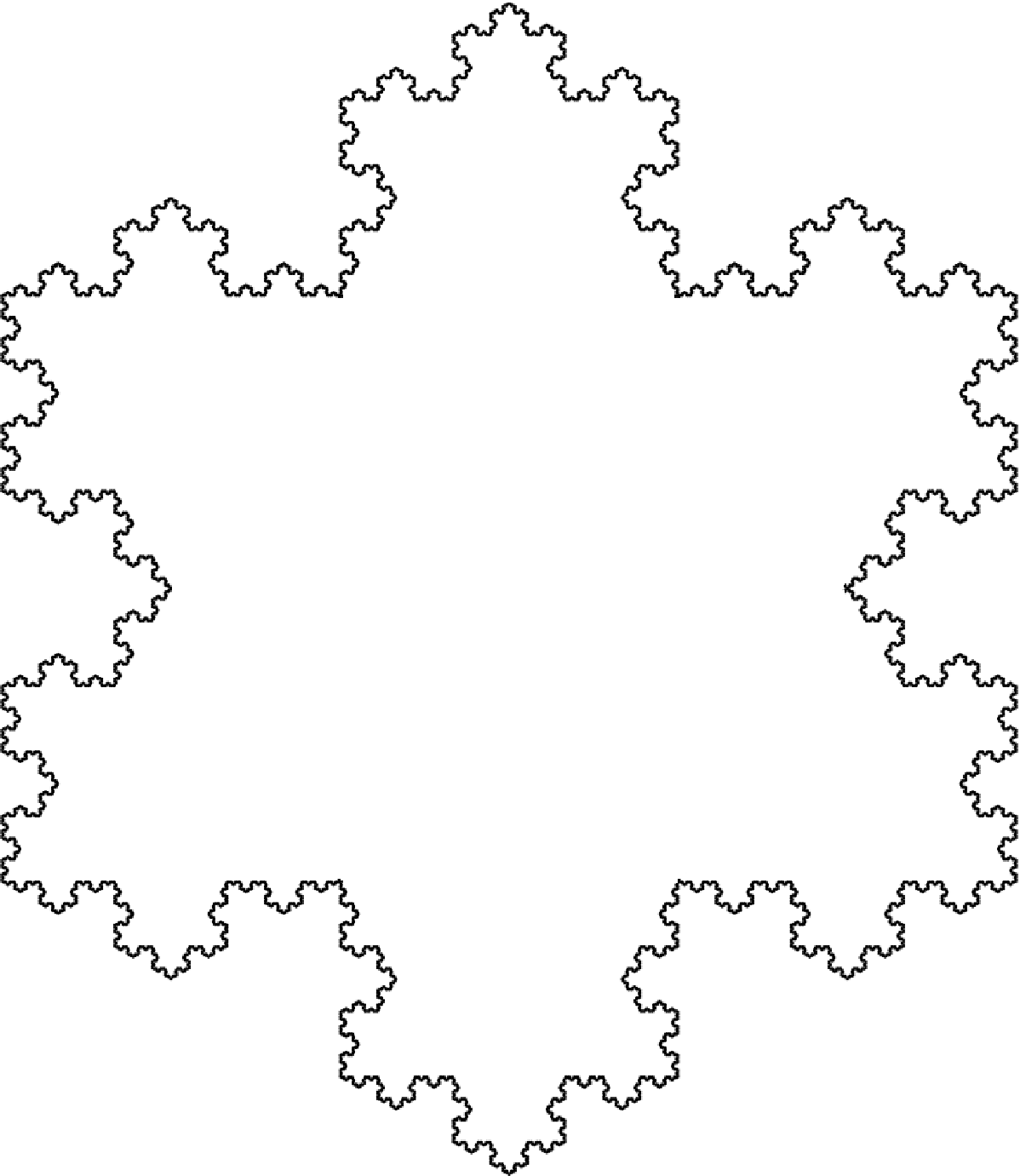}}
\hspace{5mm}
\subfloat[]{
\includegraphics[height=50mm]{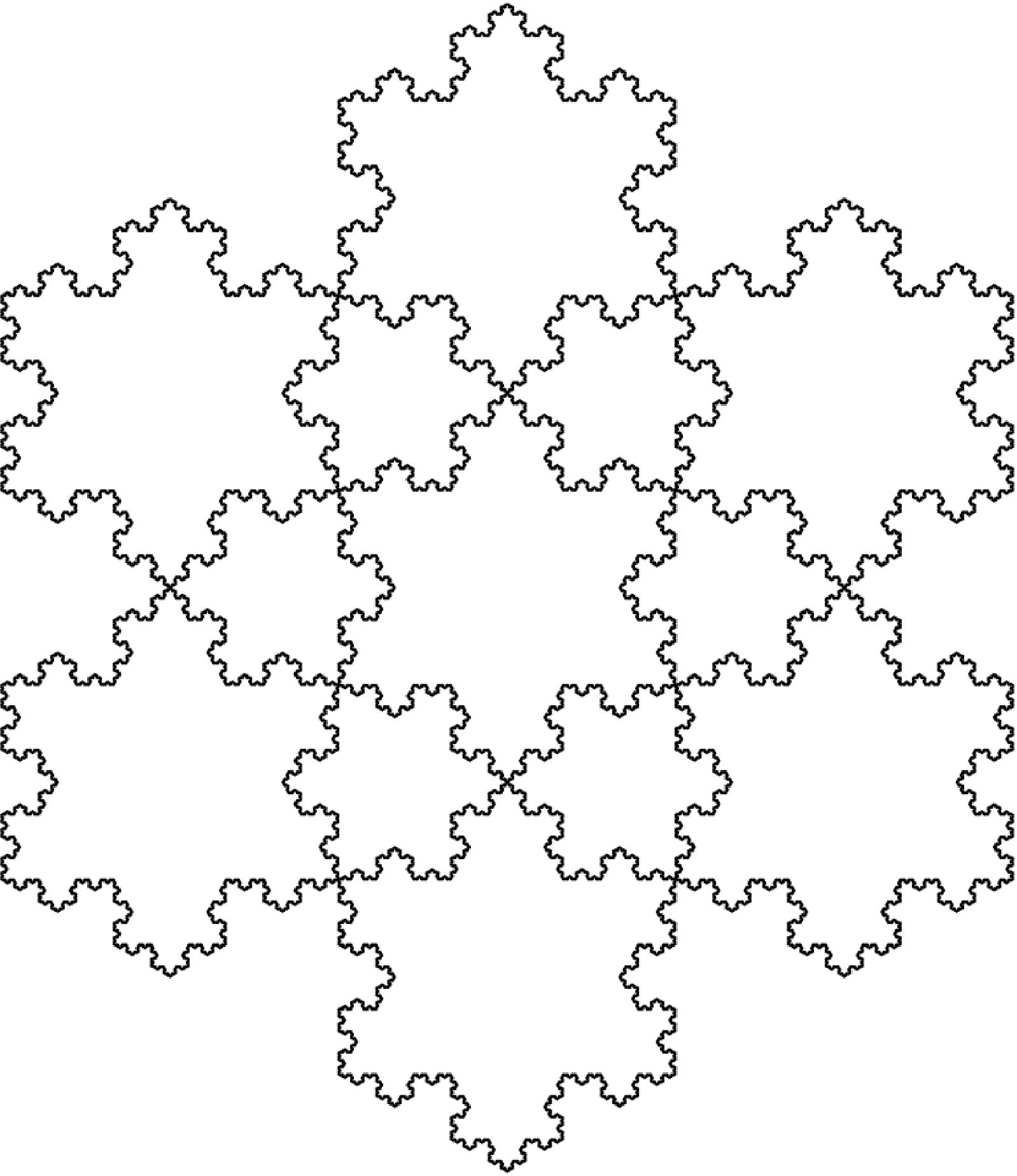}}
\hspace{5mm}
\subfloat[]{
\includegraphics[height=50mm]{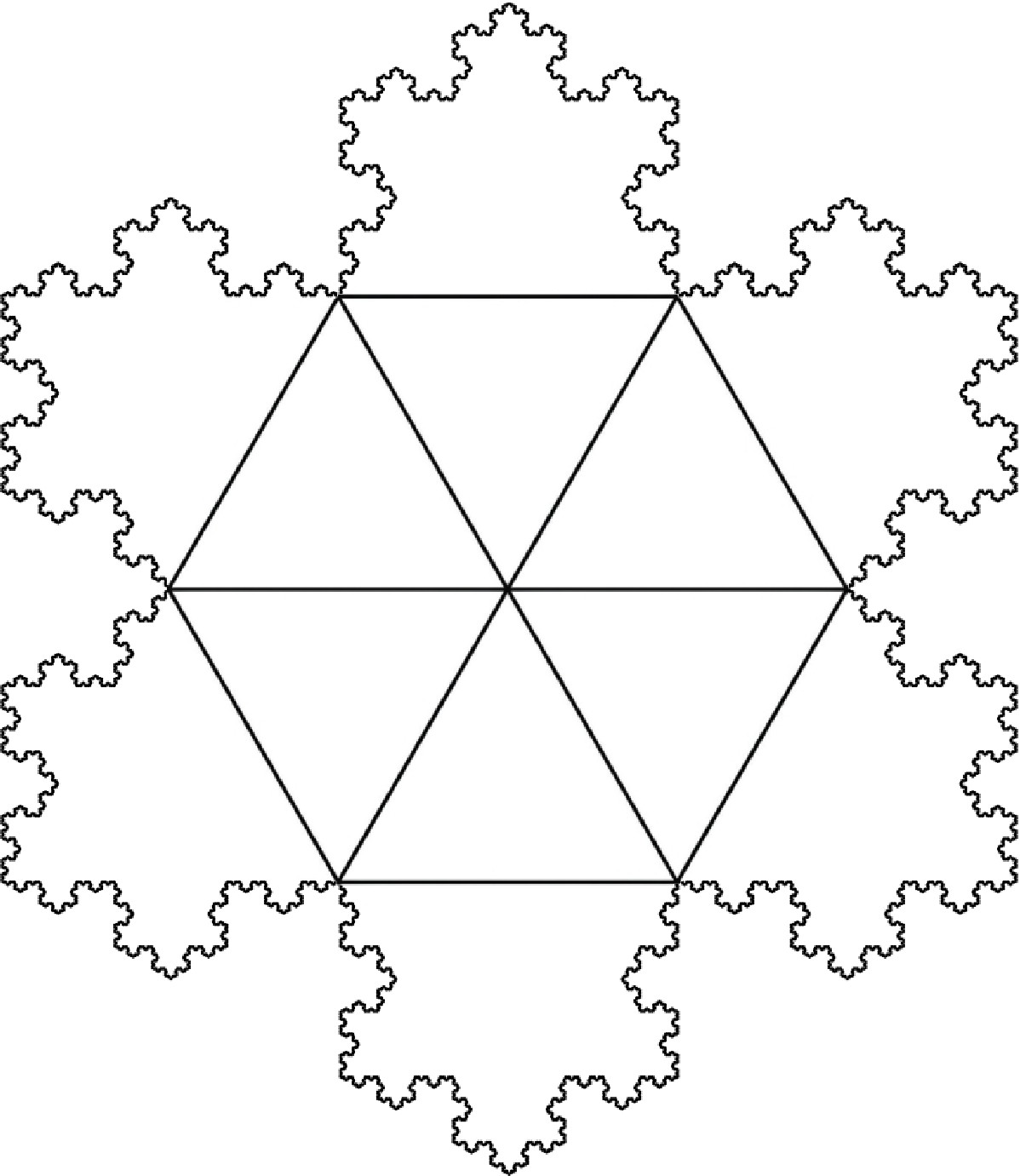}}
\caption{(a) The Koch snowflake domain $\Omega$, an open set with fractal boundary of Hausdorff dimension $\log4/\log3$. (b) A mesh of $\Omega$ comprising 
13 elements with fractal boundary (7 large and 6 small), 
each a scaled/rotated/translated version of $\Omega$ (see \cite{CaChHe25} for details). 
(c) A mesh of $\Omega$ comprising 12 elements, 6 of which are triangles and 6 of which have a boundary that is the union of a line segment and a fractal curve. 
}
\label{f:mesh}
\end{figure}

In recent years there has been increasing interest in the solution of PDEs on non-Lipschitz domains with fractal boundary \cite{bagnerini2006finite, Achdou07, capitanelli2015weighted, 
hinz2021non, creo2022transmission}, such as the Koch snowflake illustrated in Figure \ref{f:mesh}(a). Such domains cannot be meshed with a finite number of Lipschitz elements, so the results of \cite{cangiani2022hp} do not apply. An obvious remedy is to first replace the domain $\Omega$ with a smoother ``prefractal'' approximation, and mesh this. However, this can introduce significant geometry approximation errors \cite{FractalTransmission}; 
furthermore, meshing the prefractal approximation may require a large number of elements, particularly if one wants to work with shape-regular simplicial meshes. For an example of the extreme mesh complexity that can arise in such geometries with conforming FEM, see, e.g.,\ \cite{Bagnerini13}. 

To mesh a domain like that in Figure \ref{f:mesh}(a) more efficiently one might consider working with meshes containing elements with fractal boundaries, such as those illustrated in Figure \ref{f:mesh}(b),(c), so that the geometry of $\Omega$ is fully captured.
Such meshes have been successfully used in the context of integral equation methods for acoustic wave scattering by fractal structures in \cite{FractalTransmission, 
HausdorffBEM,HausdorffDomain, CaChHe25}, using the numerical quadrature rules on fractal meshes presented in \cite{HausdorffQuadrature,NonDisjointQuad}. 
However, 
aside from the results for piecewise constant approximation presented in 
\cite{HausdorffBEM,HausdorffDomain, CaChHe25}, 
piecewise polynomial approximation theory on such meshes does not seem to be available in the literature. 

In this paper we provide best approximation error estimates for discontinuous piecewise polynomial approximation on 
non-Lipschitz meshes of non-Lipschitz domains. 
In fact, many of our results hold for arbitrary meshes of arbitrary domains, meaning that 
they require  
no regularity assumption whatsoever on the boundary of $\Omega$, or on the boundaries of the mesh elements. In particular, 
these boundaries 
can be fractal {(as in Figure \ref{f:mesh})}, 
or can even have positive Lebesgue measure.

Our estimates are stated in the context of Sobolev spaces, and when dealing with non-smooth domains it is important to clarify the precise Sobolev space setting in which one is working. 
In the FEM setting one typically requires %
estimates in ``intrinsic'' Sobolev spaces such as $W^m(\Omega)$, for some $m\in\N$, defined in terms of square integrability of weak derivatives on $\Omega$. In the integral equation setting one often requires estimates in ``extrinsic'' spaces derived from the Bessel potential Sobolev space $H^s(\R^n)$, for some $s\in\R$, 
such as the space $H^s(\Omega)$ of restrictions to $\Omega$ of elements of $H^s(\R^n)$, or the space $\tH^s(\Omega)$ defined as the closure of $C^\infty_0(\Omega)$ in $H^s(\R^n)$. 
Here we provide results in both intrinsic and extrinsic spaces.

Our results are stated in \S\ref{s:results} and the proofs can be found in \S\ref{s:proofs}. Our results can be summarised thus: 
\begin{itemize}
\item
In Theorem \ref{thm:hp} we provide optimal $hp$ local best approximation error estimates based on local regularity assumptions. Theorem \ref{thm:hp} generalises \cite[Lemma 4.31, Eqn~(4.24)]{cangiani2022hp} to the non-Lipschitz case, and its proof, like that of \cite[Lemma 4.31, Eqn~(4.24)]{cangiani2022hp}, 
makes use of a ``covering mesh'' approach. However, as is explained 
{in Remark \ref{r:Cangiani}}, 
Theorem \ref{thm:hp} provides a slightly stronger result than \cite[Lemma 4.31, Eqn~(4.24)]{cangiani2022hp} {(even in the Lipschitz case)}. This strengthening means that, when using Theorem \ref{thm:hp} to derive global error estimates based on global regularity assumptions, in contrast to \cite{cangiani2022hp}, we do not need to impose any condition analogous to {the first part of} \cite[Assumption 4.28]{cangiani2022hp} {(restated as \eqref{e:Ass428a} below)}, restricting the number of mesh elements that can intersect each covering element.

\item In Theorem \ref{thm:Approx}, which is our main result, we provide 
global error estimates in broken Sobolev norms for functions $u$ in the extrinsic spaces $H^m(\Omega)$ 
for arbitrary $m\in\N$   
and for arbitrary meshes of arbitrary domains. 
\item 
In Corollary \ref{cor:Intrinsic} we provide a version of Theorem \ref{thm:Approx} valid for functions $u$ in the intrinsic space $W^m(\Omega)$. This result holds under the assumption that $\Omega$ is a $W^m$ extension domain {(see the discussion after \eqref{e:HsDef})}, which is an implicit regularity assumption on $\Omega$ (but not on the mesh). It holds whenever $\Omega$ is an $(\epsilon,\delta)$ locally uniform domain, so in particular for the Koch snowflake domain in Figure \ref{f:mesh}(a). 
{It is not guaranteed to hold if the $W^m$ extension property fails, e.g.\ if $\Omega$ has an outward power-type cusp or an interior slit.} 

\item
In Corollary \ref{cor:int} we 
derive $L^2(\Omega)$ estimates for functions in the fractional order spaces $H^s(\Omega)$ for general $s\in[0,m]$, and in Corollary \ref{cor:Duality} we provide estimates in negative order norms. These two corollaries are both valid for arbitrary meshes of arbitrary domains. 

\item Finally, in Corollary \ref{cor:Duality2} we generalise the results of Corollary \ref{cor:Duality} to distributions of negative smoothness. This result requires the additional assumption that $\{H^s(\Omega)\}_{s\in[0,m]}$ (equivalently, $\{\tH^{-s}(\Omega)\}_{s\in[0,m]}$) is an interpolation scale (in the sense of \cite[Rem.~3.8]{InterpolationCWHM}), which again is an implicit regularity assumption on $\Omega$ (but not on the mesh). Conditions under which this holds are discussed in \S\ref{s:results}. We note that it holds for any $(\epsilon,\delta)$ locally uniform domain, so in particular for the Koch snowflake domain in Figure \ref{f:mesh}(a). 
{But it also holds for a wider class of domains - see Remark \ref{rem:intscale} for examples.}
\end{itemize}

As hinted at above, 
our work is largely motivated by recent developments in the numerical solution of integral equations on fractal domains \cite{HausdorffBEM,HausdorffDomain, CaChHe25,FractalTransmission}.
Indeed, a simplified piecewise constant version of Corollary \ref{cor:Duality} was used to analyse the convergence of a piecewise constant Galerkin boundary element method for scattering by fractal screens in \cite{CaChHe25}. Furthermore, the result of Corollary \ref{cor:int} {has been recently applied} in \cite{FractalTransmission} to analyse %
piecewise polynomial volume integral equation methods for the Lippmann-Schwinger equation for scattering by fractal inhomogeneities.  

Our results are also of relevance to the analysis of dG-FEM on non-smooth domains, and in \cite{GoMoHe} we present and analyse a dG-FEM 
on 
a fractal domain 
using meshes with fractal elements of the type shown in Figure \ref{f:mesh}(b). 
We note that, for dG-FEM applications, in addition to error estimates for approximation on the mesh elements, 
one also typically needs error estimates for approximation on element boundaries (cf., e.g., \cite[Eqn (4.25)]{cangiani2022hp}). In principle, such estimates can be derived for non-Lipschitz meshes of non-Lipschitz domains, provided that 
the boundaries are sufficiently regular to ensure that suitable trace operators exist. 
This 
holds in particular for the mesh of Figure \ref{f:mesh}(b), because the Koch snowflake domain is a $W^m(\Omega)$ extension domain and its boundary is a $d$-set (in the terminology of \cite{JoWa84}) of Hausdorff dimension $d=\log{4}/\log{3}$, so by results in \cite{JoWa84} there exists a bounded linear trace operator from $W^1(\Omega)$ to $L^2(\partial\Omega,\cH^d)$, where $\cH^d$ is the $d$-dimensional Hausdorff measure. 
However, we defer further discussion of this to 
\cite{GoMoHe} and to future work.

\clearpage

\section{Definitions and results}
\label{s:results}
Let $\Omega\subset\R^n$ be a domain, {by which we mean a non-empty open set (not necessarily connected)}. For $m\in \N$ let $W^m(\Omega)$ denote the Sobolev space of functions $u\in L^2(\Omega)$ with square-integrable weak derivatives up to order $m$, normed by 
$\|u\|_{W^m(\Omega)}=\left( \sum_{j=0}^m |u|^2_{W^j(\Omega)}\right)^{1/2}$, where $|u|_{W^j(\Omega)} = \left(\sum_{|\alpha|=j}\|D^\alpha u\|^2_{L^2(\Omega)}\right)^{1/2}$.
For $s\in \R$ let $H^s(\R^n)$ denote the Bessel potential Sobolev space of tempered distributions $u$ on $\R^n$ for which 
$\|u\|_{H^s(\R^n)}:=(\int_{\R^n}(1+|\xi|^2)^s|\hat{u}(\xi)|^2\,\rd \xi)^{1/2}<\infty$. 
Let $\tH^s(\Omega)$ denote the closure of $C_0^\infty(\Omega)$ in $H^s(\R^n)$, a closed subspace of $H^s(\R^n)$, 
and let 
\begin{align}
\label{}
H^s(\Omega):=\{U|_\Omega:U\in H^s(\R^n)\}
\end{align}
{denote the space of restrictions to $\Omega$ of elements of $H^s(\R^n)$,} 
equipped with the quotient norm
\begin{align}
\label{e:HsDef}
\|u\|_{H^s(\Omega)} := \min_{\substack{U\in H^s(\R^n)\\ u=U|_\Omega}}\|U\|_{H^s(\R^n)}, \qquad u \in H^s(\Omega).
\end{align}
We recall that for any domain $\Omega\subset\R^n$ and any $s\in\R$ the dual space of $H^s(\Omega)$ can be unitarily realised as the space $\tH^{-s}(\Omega)$ (with $L^2(\Omega)$ as pivot space) \cite[Thm~3.3]{ChaHewMoi:13}. 
We also recall that, for $m\in \N_0$, %
$H^m(\R^n)=W^m(\R^n)$ with equivalent norms, and $H^m(\Omega)$ is continuously embedded in $W^m(\Omega)$. In general $H^m(\Omega)$ may be a strict subset of $W^m(\Omega)$. However, if $\Omega$ is a $W^m$ extension domain, meaning that there exists a bounded linear operator $E^m_\Omega :W^m(\Omega)\to W^m(\R^n)$ such that $E^m_\Omega (u)|_\Omega=u$ for all $u\in W^m(\Omega)$, then $H^m(\Omega)=W^m(\Omega)$ with equivalent norms. This holds in particular if $\Omega$ is Lipschitz \cite{Calderon61,Rychkov}, or more generally if $\Omega$ is $(\epsilon,\delta)$ locally uniform \cite{Jones,Rogers}. 
The latter class includes a large subclass of domains with fractal boundary; in particular, the Koch snowflake domain pictured in Figure \ref{f:mesh} is $(\epsilon,\infty)$ locally uniform \cite{Jones}, and hence is a $W^m$ extension domain. 

{We now define what we mean by a mesh. Note that $\Omega$ is not required to be bounded.}

\begin{defn}
\label{d:OmegaMesh}
Given a domain $\Omega\subset\R^n$, we say that a countable (e.g.\ finite) collection $\cT=\{K\}$ of non-empty bounded open subsets of $\Omega$ is a \emph{mesh} of $\Omega$ if $K\cap K'=\emptyset$ for $K,K'\in \cT$, $K\neq K'$, and 
\begin{align}
\label{e:Mesh}
\Big|\Omega\setminus \bigcup_{K\in \cT} K\Big|=0.
\end{align}
Given such a mesh, let $h_K:=\diam(K)$ for each $K\in \cT$, and let $h:=\sup_{K\in\cT}h_K$ denote the \emph{meshwidth} of $\cT$.
\end{defn}

As in \cite{cangiani2022hp}, we will make use of the notion of a \emph{covering} of a mesh.

\begin{defn}
Given a domain $\Omega\subset\R^n$, and a mesh $\cT=\{K\}$ of $\Omega$, we define a \emph{covering} $\cT^\#=\{\cK\}$ of $\cT$ to be a set of open simplices or {parallelotopes} $\cK$ such that for each $K\in \cT$ there exists at least one $\cK\in\cT^\#$ with $K\subset\cK$. Given a covering $\cT^\#$, {for each $\cK\in\cT^\#$ we define $h_\cK:=\diam(\cK)$ and $\sigma_{\cK}:=h_\cK/\rho_\cK$, where $\rho_\cK$ is the diameter of the largest ball that can be inscribed in $\cK$.} 
We say that $\kappa:\cT\to \cT^\#$ is a \emph{covering choice function} if $K\subset\kappa(K)$ for each $K\in \cT$. 
\end{defn}

{In Theorem \ref{thm:hp}, which is our first main result, we provide} 
an optimal local $hp$ approximation result, generalising \cite[Lemma 4.31, Eqn~(4.24)]{cangiani2022hp} to the non-Lipschitz case. 
{A detailed comparison between Theorem \ref{thm:hp} and \cite[Lemma 4.31, Eqn~(4.24)]{cangiani2022hp} is provided in Remark \ref{r:Cangiani}.} 
{We emphasize that Theorem \ref{thm:hp} holds for arbitrary meshes of arbitrary domains.} 

Here, and henceforth, given $p\in\N_0$ and $S\subset\R^n$ let $\cP_p(S)$ denote the set of polynomials on $S$ in $n$ variables of total degree less than or equal to $p$. 

\begin{thm}
\label{thm:hp}
Let $\Omega\subset\R^n$ be a domain, $\cT$ a mesh of $\Omega$ with {meshwidth $h>0$, $\cT^\#$ a covering of $\cT$ for which $h\leq h_0:=\sup_{\cK\in\cT^\#}\diam(\cK)<\infty$,}  and $\kappa:\cT\to \cT^\#$ a covering choice function. 
Let $u\in L^2(\Omega)$ and suppose that for each $\cK\in\cT^\#$ there exists $m_\cK\in\N_0$ and $U_\cK\in W^{m_\cK}
(\cK)$ such that $(U_\cK)|_K=u|_K$ for each 
$K\in \kappa^{-1}(\cK)$. Then for each $\cK\in\cT^\#$ and each $p\in\N_0$ there exists an operator $\pi_{p,m_\cK,\cK}:W^{m_\cK}(\cK)\to \cP_p(\cK)$ and a constant $C_{n,h_0,m_\cK,{\sigma_\cK}}>0$, depending only on $n$, $h_0$, $m_\cK$ and 
{$\sigma_\cK$}, 
such that, with $t_\cK:=\min(m_\cK,p+1)$, 
\begin{align}
\label{e:Est1}
\Bigg(\sum_{K\in \kappa^{-1}(\cK)} \|u-\pi_{p,m_\cK,\cK}(U_\cK)\|^2_{W^j(K)}\Bigg)^{1/2}
\leq C_{n,h_0,m_\cK,{\sigma_\cK}} \frac{h_\cK^{t_\cK-{j}}}{(p+1)^{m_\cK-j}} \sum_{r=t_\cK}^{m_\cK} 
|U_\cK|_{W^r({\cK})},
\quad j=0,1,\ldots,m_\cK.
\end{align}
\end{thm}

{
\begin{rem}
\label{r:Cangiani}
Here 
we compare Theorem \ref{thm:hp} with the related result presented in \cite[Lemma 4.31]{cangiani2022hp}.

We begin by re-emphasizing that Theorem \ref{thm:hp} holds for arbitrary meshes of arbitrary domains, whereas \cite[Lemma 4.31]{cangiani2022hp} was stated only for a class of Lipschitz meshes of Lipschitz domains. 

We note further that 
our result 
strengthens the statement of \cite[Lemma 4.31, Eqn~(4.24)]{cangiani2022hp} (including in the Lipschitz case), because in \eqref{e:Est1} we control in one estimate the contributions from all the elements $K$ that are associated, via the {covering choice} function $\kappa$, with a given covering element $\cK$, rather than providing a separate bound for each such $K$ as in \cite[Lemma 4.31, Eqn~(4.24)]{cangiani2022hp}. 
This means that, when deriving the global approximation error estimates presented in Theorem \ref{thm:Approx} below, we avoid the need for a mesh density assumption like the first part of \cite[Assumption 4.28]{cangiani2022hp},  
viz.
\begin{align}
\label{e:Ass428a}
\exists\, \mathcal{O}_\Omega\in \N, \,\,
\max_{K\in \cT} {\rm card} \{ K'\in \cT:K'\cap \cK\neq \emptyset, 
\,\,\cK\in \cT^\# \text{ such that } K\subset\cK  \} \leq \mathcal{O}_\Omega,
\end{align}
which was included in the hypothesis of \cite[Lemma 4.31]{cangiani2022hp}. (The condition \eqref{e:Ass428a} is not used in the proof of \cite[Lemma 4.31, Eqn~(4.24)]{cangiani2022hp} but is implicitly used later in \cite{cangiani2022hp} in the derivation of global error estimates based on \cite[Lemma 4.31]{cangiani2022hp}, see e.g.\ \cite[Thm~5.6]{cangiani2022hp} and the remarks immediately following it.) 

The results we derive in Theorem \ref{thm:Approx} (without assuming \eqref{e:Ass428a}) are suited to the analysis of numerical methods based on either $h$ or $p$ refinement, in cases where one has a global regularity assumption on the function $u$ being approximated. 
For the analysis of $hp$ methods on locally refined meshes, as in \cite{cangiani2022hp}, it is preferable to have element-wise estimates that involve local regularity assumptions on $u$. Such estimates can be trivially derived from our Theorem \ref{thm:hp} by simply bounding below the sum on the left-hand side of \eqref{e:Est1} by any one of its summands. 
In this context it is natural to restrict attention to covering meshes $\cT^\#$ that are well adapted to the local refinement of the mesh $\cT$. This can be achieved by imposing an assumption like the second part of \cite[Assumption 4.28]{cangiani2022hp} 
(which was included in the hypothesis of \cite[Lemma 4.31]{cangiani2022hp}), 
viz.
\begin{align}
\label{e:Ass428b}
\exists\, C_{\rm diam}>0, \,\,
h_\cK\leq C_{\rm diam} h_K \text{ for all } K\in \cT \text{ and } \cK\in \cT^\#  \text{ such that } K\subset\cK.
\end{align}
The point of assumption \eqref{e:Ass428b} is that it allows one to replace the covering element meshwidth $h_\cK$ on the right-hand side of \eqref{e:Est1} by the element meshwidth $h_K$, 
to obtain that, for any $K\in \kappa^{-1}(\cK)$,
\begin{align}
\label{e:Est1single}
\|u-\pi_{p,m_\cK,\cK}(U_\cK)\|_{W^j(K)}
\leq C'_{n,h_0,m_\cK,\sigma_\cK} \frac{h_K^{t_\cK-j}}{(p+1)^{m_\cK-j}} \sum_{r=t_\cK}^{m_\cK} 
|U_\cK|_{W^r(\cK)},
\quad j=0,1,\ldots,m_\cK,
\end{align}
where $C'_{n,h_0,m_\cK,\sigma_\cK}$ may differ from $C_{n,h_0,m_\cK,\sigma_\cK}$ as a result of the application of the bound in \eqref{e:Ass428b}. 

To show that Theorem \ref{thm:hp} implies the result presented in \cite[Lemma 4.31, Eqn~(4.24)]{cangiani2022hp}, 
we note that 
Theorem \ref{thm:hp} assumes the existence of $U_\cK\in W^{m_\cK}
(\cK)$ such that $(U_\cK)|_K=u|_K$. When 
$\Omega$ is 
Lipschitz, there exists a universal extension operator $E_\Omega$ which maps $W^{m}(\Omega)\to W^{m}(\R^n)$ continuously for all $m\in\N_0$ \cite{Calderon61,Stein}. 
In this case, assuming that $E_\Omega(u)|_{\cK}\in W^{m_\cK}
(\cK)$, one can take 
$U_\cK=E_\Omega(u)|_{\cK}$ in \eqref{e:Est1single} and, bounding the sum of semi-norms on the right-hand side of \eqref{e:Est1single} by the $W^{m_\cK}(\Omega)$ norm, we obtain
\begin{align}
\label{e:Est1_cang}
 \|u-\pi_{p,m_\cK,\cK}(E_\Omega (u))\|_{W^j(K)}
\leq C'_{n,h_0,m_\cK,\sigma_\cK} \frac{h_K^{t_\cK-j}}{(p+1)^{m_\cK-j}} \|E_\Omega (u)\|_{W^{m_\cK}(\cK)},
\quad j=0,1,\ldots,m_\cK,
\end{align}
which recovers the result presented in \cite[Lemma 4.31, Eqn~(4.24)]{cangiani2022hp}.\footnote{{Equation  \eqref{e:Est1_cang} corrects a typo in \cite[Lemma 4.31, Eqn~(4.24)]{cangiani2022hp}, where we believe that ``$\pi_p v$'' should be replaced by ``$\pi_p(\mathfrak{E}v)$'' (in the notation of \cite{cangiani2022hp}).}}  %

One further difference between our result and that of \cite[Lemma 4.31, Eqn~(4.24)]{cangiani2022hp} 
relates to the dependencies of the constants involved. 
The constant $C_{n,h_0,m_\cK,\sigma_\cK}$ in \eqref{e:Est1} depends only on 
$n$, $m_\cK$, $\sigma_\cK$ and $h_0=\sup_{\cK\in\cT^\#}h_\cK$, while the corresponding constant $C_1$ in 
\cite[Lemma 4.31, Eqn~(4.24)]{cangiani2022hp} depends on $n$, $m_\cK$, $\sigma_\cK$, $j$, $C_{\rm diam}$, and on $\Omega$. 
The dependence of our constant $C_{n,h_0,m_\cK,\sigma_\cK}$ on $h_0$ arises because the standard $hp$ polynomial approximation theory results used in the proof of both Theorem \ref{thm:hp} and \cite[Lemma 4.31, Eqn~(4.24)]{cangiani2022hp} are sharp in their $h$ dependence in the limit as $h\to 0$, but not in the limit as $h\to \infty$. So the constant in our estimate can only be independent of $h_\cK$ if $h_0$ is finite. If one assumes $\Omega$ is bounded, as in \cite{cangiani2022hp}, then the meshwidth of $\cT$ automatically satisfies the upper bound $h\leq \diam(\Omega)$. If one assumes further that \eqref{e:Ass428b} holds, then $h_0$ is necessarily finite, with $h_0\leq C_{\rm diam}\diam(\Omega)$. Combined with the above observations, this clarifies that the dependence of the constant $C_1$ in \cite[Lemma 4.31, Eqn~(4.24)]{cangiani2022hp} on $\Omega$ can be assumed to be through $\diam(\Omega)$ alone, not on any other properties of $\Omega$. Our need to explicitly assume $h_0$ is finite in Theorem \ref{thm:hp} arises because we do not assume that $\Omega$ is bounded, nor that \eqref{e:Ass428b} holds. 
\end{rem}
}

Theorem \ref{thm:hp} is the foundation of 
our {second} main result, Theorem \ref{thm:Approx}, which provides global best approximation error estimates in broken Sobolev norms for $h$ and $p$ refinements for functions in the extrinsic spaces $H^m(\Omega)$ {(which are guaranteed, by the definition of $H^m(\Omega)$, to be extendable to $\R^n$)}. 
We emphasize that Theorem \ref{thm:Approx} holds for arbitrary meshes of arbitrary domains.

\begin{defn}
Given a domain $\Omega\subset\R^n$, a mesh $\cT=\{K\}$ of $\Omega$ of meshwidth $h>0$, and $p\in \N_0$, 
we define the \emph{piecewise polynomial space} 
\begin{align}
\label{}
V_{h,p}:=\{u\in L^2(\Omega):u|_{K}\in \cP_{p} \text{ for each }K\in\cT\}\subset L^2(\Omega).
\end{align}

\end{defn}

\begin{thm}
\label{thm:Approx}
Let $n,m\in \N$ and $h_0>0$. There exists a constant $C_{n,h_0,m}>0$, depending only on $n$, $h_0$ and $m$, such that if $\Omega\subset\R^n$ is a domain, $\cT$ is a mesh of $\Omega$ with meshwidth $h\in (0,h_0]$, $u\in H^m(\Omega)$, $U\in H^m(\R^n)=W^m(\R^n)$ is such that $U|_\Omega = u$, 
and $p\in \N_0$, then, 
with $t:=\min(m,p+1)$, 
 \begin{align}
\label{e:Est1ex}
\min_{v\in V_{h,p}}\Big(\sum_{K\in\cT} \|u-v\|^2_{W^j(K)}\Big)^{1/2}
\leq C_{n,h_0,m} \frac{h^{t-j}}{(p+1)^{m-j}} \sum_{r=t}^m |U|_{W^r(\R^n)},
\qquad j=0,1,\ldots,m.
\end{align}
\end{thm}

Theorem \ref{thm:Approx} immediately implies the following corollary for functions in  the intrinsic spaces $W^m(\Omega)$, under the assumption that $\Omega\subset\R^n$ is a {$W^m$} extension domain (cf.\  the discussion after \eqref{e:HsDef}). 
\begin{cor}
\label{cor:Intrinsic}
Let $n,m\in \N$ and $h_0>0$. 
If $\Omega\subset\R^n$ is a $W^m$ extension domain and $\cT$ is a mesh of $\Omega$ with meshwidth $h\in (0,h_0]$, then, 
for each $p\in \N_0$, and with $t:=\min(m,p+1)$, 
\begin{align}
\label{e:Est3a}
\min_{v\in V_{h,p}}
\Big(\sum_{K\in\cT} \|u-v\|^2_{W^j(K)}\Big)^{1/2}
\leq C_{m,\Omega}\, C_{n,h_0,m} \frac{h^{t-j}}{(p+1)^{m-j}} \|u\|_{W^m(\Omega)},
\qquad j=0,1,\ldots,m, \quad u\in W^m(\Omega),
\end{align}
where 
$C_{n,h_0,m}$ is the constant from Theorem \ref{thm:Approx} and 
$C_{m,\Omega}$ is the norm of 
$E^m_\Omega:W^m(\Omega)\to W^m(\R^n)$.
\end{cor}

From Theorem \ref{thm:Approx} one can derive three further corollaries, motivated by the integral equation applications in \cite{HausdorffBEM,HausdorffDomain,CaChHe25,FractalTransmission}. 
The first
concerns $L^2$ approximation errors for functions in fractional order spaces. 
We mention here that while the proof of 
Corollary \ref{cor:int} (presented in \S\ref{s:proofs}) uses function space interpolation, 
it does not require the spaces $\{H^s(\Omega)\}_{s\in[0,m]}$ to be an interpolation scale;\footnote{This is known not to hold for arbitrary $\Omega$ - see, e.g., \cite[Lem.~4.8]{InterpolationCWHM}.} it requires only that $H^{\theta m}(\Omega)$ is continuously embedded in the interpolation space $K_{\theta,2}((L^2(\Omega),H^m(\Omega)))$ for every $\theta\in(0,1)$, which holds for all open $\Omega$ by \cite[Lem.~4.2]{InterpolationCWHM}.
\begin{cor}
\label{cor:int}
Let $n,m\in \N$ and $h_0>0$. 
Suppose that $\Omega\subset\R^n$ is a domain, $\cT$ is a mesh of $\Omega$ with meshwidth $h\in (0,h_0]$, $p\geq m-1$, 
and $0\leq s\leq m$ (with $s$ not necessarily an integer). Then 
\begin{align}
\label{e:Est2}
\min_{v\in V_{h,p}}\|u-v\|_{L^2(\Omega)} 
\leq (C_{n,h_0,m})^{s/m}\left(\frac{\,h}{p+1}\right)^s \|u\|_{H^s(\Omega)}
,
\qquad u\in H^s(\Omega).
\end{align} 
\end{cor}

The second, presented as Corollary \ref{cor:Duality} below,
provides best approximation error estimates in negative order norms. 
To aid the interpretation 
of Corollary \ref{cor:Duality}, 
we recall that, for any domain $\Omega\subset\R^n$, extension by zero is a unitary isomorphism from $H^0(\Omega)=L^2(\Omega)$ to $\tH^0(\Omega)$ (its inverse being the restriction operator). The fact that $\tH^{s_1}(\Omega)\supset \tH^{s_2}(\Omega)$ for $s_1<s_2$ then implies that extension by zero provides a continuous embedding of $H^0(\Omega)=L^2(\Omega)$ into $\tH^{s}(\Omega)$ for $s\leq 0$, a fact we are using implicitly in 
what follows. 

\begin{cor}
\label{cor:Duality}
Let $n,m\in \N$ and $h_0>0$. 
Suppose that $\Omega\subset\R^n$ is a domain, $\cT$ is a mesh of $\Omega$ with meshwidth $h\in (0,h_0]$, $p\geq m-1$, 
and $-m\leq s_1\leq 0 \leq s_2\leq m$ (with $s_1,s_2$ not necessarily integers). Then 
\begin{align}
\label{eq:Hsboundh_dn}
\min_{v\in V_{h,p}}
\|u-v\|_{\tH^{s_1}(\Omega)}
\leq (C_{n,h_0,m})^{(s_2-s_1)/m}\left(\frac{h}{p+1}\right)^{s_2-s_1}\|u\|_{H^{s_2}(\Omega)},\qquad u\in H^{s_2}(\Omega).
\end{align}

\end{cor}

The third extends Corollary \ref{cor:Duality} to functions/distributions in negative order spaces. 

\begin{cor}
\label{cor:Duality2}
Let $n,m\in \N$ and $h_0>0$. 
Suppose that $\Omega\subset\R^n$ is a domain, and 
that $\{H^s(\Omega)\}_{0\leq s\leq m}$ is an interpolation scale (equivalently, that $\{\tH^s(\Omega)\}_{-m\leq s\leq 0}$ is an interpolation scale). 
Suppose that $\cT$ is a mesh of $\Omega$ with meshwidth $h\in (0,h_0]$, $p\geq m-1$, and that $-m\leq s_1\leq s_2\leq 0$ (with $s_1$ and $s_2$ not necessarily integers). 
Then 
\begin{align}
\label{eq:Hsboundh_dn2}
\min_{v\in V_{h,p}}
\|u-v\|_{\tH^{s_1}(\Omega)}\leq \tilde{C}(C_{n,h_0,m})^{(s_2-s_1)/m}\left(\frac{h}{p+1}\right)^{s_2-s_1}\|u\|_{\tH^{s_2}(\Omega)},\qquad u\in \tH^{s_2}(\Omega),
\end{align}
where $\tilde{C}$ is 
the norm\footnote{This norm is not guraranteed to equal one, since even when $\{H^s(\Omega)\}_{0\leq s\leq m}$ is an interpolation scale, it is not in general exact  - see, e.g., \cite[Lem.~4.13]{InterpolationCWHM}.} of the embedding of $\tH^{s_2}(\Omega)$ into the interpolation space $K_{1-s_2/s_1,2}(\tH^{s_1}(\Omega),\tH^0(\Omega))$.
\end{cor}

\begin{rem}
\label{rem:intscale}
Regarding the interpolation scale assumption in  Corollary \ref{cor:Duality2}:
\begin{enumerate}[(i)]
\item We recall from \cite[Cor.~4.9]{InterpolationCWHM} 
that $\{H^s(\Omega)\}_{0\leq s\leq m}$ is an interpolation scale (in the sense of \cite[Rem.~3.8]{InterpolationCWHM}) if and only if $\{\tH^s(\Omega)\}_{-m\leq s\leq 0}$ is an interpolation scale. 
\item 
If $\Omega\subset\R^n$ is a Lipschitz hypograph or a Lipschitz domain then $\{H^{s}(\Omega)\}_{s\in\R}$ and $\{\tH^{s}(\Omega)\}_{s\in\R}$ are interpolation scales (see, e.g., \cite[Corollaries 4.7 and 4.10]{InterpolationCWHM}).
\item If $\Omega$ is an $(\epsilon,\delta)$ locally uniform domain then $\{H^{s}(\Omega)\}_{s\geq 0}$ and $\{\tH^{s}(\Omega)\}_{s\leq 0}$ are interpolation scales (see, e.g., \cite[Corollaries 4.7 and 4.10]{InterpolationCWHM}). 
In particular, this includes the Koch snowflake domain pictured in Figure \ref{f:mesh}. 
\item 
It was shown recently in \cite[Cor.~3.20]{CaChHe25} that $\{H^{s}(\Omega)\}_{s\in\R}$ and $\{\tH^{s}(\Omega)\}_{s\in\R}$ are interpolation scales provided $\Omega$ is ``thick'' in the sense of Triebel (see \cite[Def.~C.3]{CaChHe25}), $\overline\Omega\neq \R^n$ and the Assouad dimension of $\partial\Omega$ (as defined in e.g.\ \cite[Appendix F]{CaChHe25}) is less than $n$. 
In particular, these hypotheses hold if $\Omega$ is the interior of an $n$-attractor (in the terminology of \cite{CaChHe25}), i.e.\ the attractor of an iterated function system of contracting similarities satisfying the open set condition and having Hausdorff dimension $n$. 
This applies to all the examples in \cite[Fig.~1]{CaChHe25} (including the Koch snowflake, the Gosper island and the Levy dragon). 
These hypotheses also hold for all the ``classical snowflakes'' in \cite[\S5.1]{caetano2019density} (which includes the Koch snowflake). 
\item By combining \cite[Prop.~3.15]{CaChHe25} with \cite[Rem.~3.23]{CaChHe25} and the Wolff interpolation Theorem \cite[Thm~1]{Wolff1982}, one can show that $\{H^{s}(\Omega)\}_{s\geq 0}$ and $\{\tH^{s}(\Omega)\}_{s\leq 0}$ are interpolation scales whenever $\Omega$ is an open $n$-set (in the terminology of \cite{CaChHe25}, called an \textit{interior regular domain} in \cite{caetano2019density}), $\overline\Omega\neq \R^n$ and the Assouad dimension of $\partial\Omega$ is less than $n$. 
\end{enumerate}
Some further references to the literature on interpolation scale results for $H^s(\Omega)$ and $\tilde{H}^s(\Omega)$, and more general Besov and Triebel-Lizorkin spaces, can be found in \cite{CaChHe25}.

\end{rem}

\section{Proofs}

\label{s:proofs}
\subsection*{Proof of Theorem \ref{thm:hp}}

Our proof is similar to that of \cite[Lemma 4.31]{cangiani2022hp}, except that we prove a slightly stronger result by collecting together the contributions of all the mesh elements associated with a given covering element. 
Under the hypotheses of the lemma, for each $\cK\in\cT^\#$ let $\pi_{p,m_\cK,\cK}:W^{m_\cK}(\cK)\to \cP_p(\cK)$ denote a known optimal $hp$-version approximation operator on the simplex/{parallelotope} $\cK$ (e.g.\ \cite{babuska1987optimal,babuska1987hp,Sc:98}, \cite[Appendix C]{smears2015})\footnote{{We mention \cite[Appendix C]{smears2015} because it highlights and rectifies the fact that the Babu\v{s}ka–Suri projector defined in \cite{babuska1987optimal,babuska1987hp} does not preserve polynomials, which is required for the derivation of optimal $hp$ bounds.}}.  
Since $(U_\cK)|_K=u$ for all $K\in \kappa^{-1}(\cK)$ and since $K\cap K'=\emptyset$ for $K\neq K'$, for $j=0,1,\ldots,m_\cK$ we can apply standard $hp$ approximation results (e.g.\ \cite{babuska1987optimal,babuska1987hp,Sc:98}, \cite[Appendix C]{smears2015}) to estimate, with $t_\cK:=\min(m_\cK,p+1)$,
\begin{align*}
\Bigg(\sum_{K\in \kappa^{-1}(\cK)} \|u-\pi_{p,m_\cK,\cK}(U_\cK)\|^2_{W^j(K)}\Bigg)^{1/2}
&=\Bigg(\sum_{K\in \kappa^{-1}(\cK)} \|U_\cK-\pi_{p,m_\cK,\cK}(U_\cK)\|^2_{W^j(K)}\Bigg)^{1/2}\\
&\leq \|{U_\cK}-\pi_{p,m_\cK,\cK}(U_\cK)\|_{W^j(\cK)}\\
&\leq C \frac{h_\cK^{t_\cK-{j}}}{(p+1)^{m_\cK-j}} \sum_{r=t_\cK}^{m_\cK} 
|U_\cK|_{W^r({\cK})},
\end{align*}
for some constant $C>0$ depending only on $n$, $m_\cK$, $h_0$ and the shape regularity {$\sigma_\cK$} of $\cK$. 

\subsection*{Proof of Theorem \ref{thm:Approx}}

Our approach is a generalisation of the argument used to prove \cite[Proposition 3.25]{CaChHe25}, which considered 
only the case $m=1$, $p=0$. The idea is to apply Theorem \ref{thm:hp}, using a covering formed of {congruent} overlapping cubes. 

Under the hypotheses of the theorem, 
for each $\bz\in \Z^n$ let $Q_{\bz}$ be the cube of side length $h$ defined by $Q_{\bz}:=[0,h]^n + h\bz$, and let $Q_{\bz}':=[-h,2h]^n+h\bz$ be the cube of side length $3h$ formed by taking the union of $Q_{\bz}$ with all its neighbouring cubes, which satisfies $\diam(Q_\bz')=3n^{1/2}h$.  
For each $K\in\cT$, since $\diam(K)\leq h$, $K$ intersects the cube $Q_\bz$ for at least one and at most $2^n$ different $\bz\in \Z^n$. And for each such $\bz$ we have $K\subset Q_{\bz}'$. 
Hence $\cT^\#:=\{Q_\bz'\}_{z\in \Z^n}$ is a covering of $\cT$. Let $\kappa:\cT \to \cT^\#$ be any covering choice function, and let $\bz_\kappa:\cT\to \Z^n$ be defined by $\kappa(K)=Q_{\bz_\kappa(K)}'$ for $K\in\cT$. 

For each $K\in \cT$ define $v|_{K}$ to be the restriction to $K$ of $\pi_{p,m,Q_{\bz_\kappa(K)}'}(U|_{Q_{\bz_\kappa(K)}'})$, where $\pi_{p,m,Q_{\bz_\kappa(K)}'}$ is the operator from Theorem \ref{thm:hp}. 
This defines 
$v$ almost everywhere on $\Omega$, because of the fact that $\left|\Omega\setminus \bigcup_{K\in \cT} K\right|=0$. 
Furthermore, {noting that since the cubes making up $\cT^\#$ are all congruent, and hence all have the same shape regularity,} by Theorem \ref{thm:hp} we have that
\begin{align*}
 \label{}
\sum_{K\in\cT}\|u - v\|^2_{W^j(K)}
 = \sum_{\bz\in \Z^n} \sum_{K\in \bz_\kappa^{-1}(\bz)}
 \|u - v\|^2_{W^j(K)}
 &\leq\sum_{\bz\in \Z^n} \left({\tilde{C}_{n,h_0,m}}\frac{(3n^{1/2}h)^{t-j}}{(p+1)^{m-j}}\right)^2 \sum_{r=t}^m |U|^2_{W^r(Q_{\bz}')}.
 \end{align*}
{for some constant} 
$\tilde{C}_{n,h_0,m}>0$ depending only on $n$, $h_0$ and $m$. Then 
\begin{align*}
\label{}
\sum_{K\in\cT}\|u - v\|^2_{W^j(K)} 
&\leq \left(\tilde{C}_{n,h_0,m}\frac{(3n^{1/2}h)^{t-j}}{(p+1)^{m-j}}\right)^2 \sum_{\bz\in \Z^n}\sum_{r=t}^m |U|^2_{W^r(Q_{\bz}')}\\
& = 3^n \left(\tilde{C}_{n,h_0,m}\frac{(3n^{1/2}h)^{t-j}}{(p+1)^{m-j}}\right)^2\sum_{\bz\in \Z^n} \sum_{r=t}^m |U|^2_{W^r(Q_{\bz})}\\
&=
3^n \left(\tilde{C}_{n,h_0,m}\frac{(3n^{1/2}h)^{t-j}}{(p+1)^{m-j}}\right)^2\sum_{r=t}^m |U|^2_{W^r(\R^n)},
\end{align*}
where we used the fact that for each $\bz$ the 
contribution $|U|^2_{W^r(Q_{\bz})}$ gets counted $3^n$ times in the sum $\sum_{\bz\in \Z^n} |U|^2_{W^r(Q_{\bz}')}$. 
From this, \eqref{e:Est1ex} follows by taking square roots, with $C_{n,h_0,m}=3^{n/2}(3n^{1/2})^m\tilde{C}_{n,h_0,m}$. 
In the case where $\Omega$ is bounded, the fact that $v\in L^2(\Omega)$, and hence that $v\in V_{h,p}$, is immediate from the definition of $v$. In the case where $\Omega$ is unbounded the fact that $\|v\|_{L^2(\Omega)}<\infty$ follows from \eqref{e:Est1ex} with $j=0$, the fact that $u\in L^2(\Omega)$, and the triangle inequality $\|v\|_{L^2(\Omega)}\leq \|u-v\|_{L^2(\Omega)} + \|u\|_{L^2(\Omega)}$.

\subsection*{Proof of Corollary \ref{cor:Intrinsic}}
Just apply \eqref{e:Est1ex} with $U=E^m_\Omega u$ and bound $\|E^m_\Omega  w\|_{W^m(\R^n)}\leq C_{m,\Omega}\|w\|_{W^m(\Omega)}$.

\subsection*{Proof of Corollary \ref{cor:int}}

Let $\Pi:L^2(\Omega)\to V_{h,p}$ denote the $L^2$-orthogonal projection onto $V_{h,p}$, which satisfies $\|u-\Pi u\|_{L^2(\Omega)}=\min_{v\in V_{h,p}}\|u-v\|_{L^2(\Omega)}$ for $u\in L^2(\Omega)$.  
For $u\in H^m(\Omega)$, taking $j=0$ in \eqref{e:Est1ex} and bounding 
$|U|_{W^m(\R^n)}\leq \|U\|_{H^m(\R^n)}$ gives 
\begin{align}
\label{e:Est1a}
\|u-\Pi u\|_{L^2(\Omega)} \leq 
C_{n,h_0,m}\left(\frac{h}{p+1}\right)^m \|U\|_{H^m(\R^n)}.
\end{align} 
Since \eqref{e:Est1a} holds for every $U\in H^m(\R^n)$ with $U|_\Omega=u$, by the definition of $\|\cdot\|_{H^m(\Omega)}$ we have
\begin{align}
\label{e:Est1b}
\|u-\Pi u\|_{L^2(\Omega)} \leq 
C_{n,h_0,m}\left(\frac{h}{p+1}\right)^m \|u\|_{H^m(\Omega)}.
\end{align} 
We note also for $u\in L^2(\Omega)=H^0(\Omega)$ the trivial estimate
\begin{align}
\label{e:Est1c}
\|u-\Pi u\|_{L^2(\Omega)} \leq 
\|u\|_{L^2(\Omega)}.
\end{align} 
We then deduce \eqref{e:Est2} by interpolation applied to the operator $I-\Pi $, using  \eqref{e:Est1b} and \eqref{e:Est1c}.  We adapt the 
argument of \cite[Eqn~(46)]{BEMfract}, where only the special case $m=1$ was considered. 
Precisely, given $s\in[0,m]$, let $\theta:=s/m\in[0,1]$ and let $\|\cdot\|_\theta$ denote the norm on the $K$-method interpolation space $K_{\theta,2}((L^2(\Omega),H^m(\Omega)))$, with the normalisation defined in \cite[(7,8)]{InterpolationCWHM}. Then, with $\sigma:=(C_{n,h_0,m})^{1/m}h/(p+1)$, 
 \begin{equation} \label{e:Interp}
\|u-\Pi  u\|_{L^2(\Omega)} \leq \sigma^{\theta m} \|u\|_\theta \leq \sigma^{\theta m} \|u\|_{H^{\theta m}(\Omega)} = \sigma^{s} \|u\|_{H^s(\Omega)},
\end{equation}
where the first inequality follows by interpolation since $K_{\theta,2}((L^2(\Omega),L^2(\Omega)))=L^2(\Omega)$ with equality of norms given the normalisation \cite[(7,8)]{InterpolationCWHM} (see \cite[Lem.~2.1(iii), Thm.~2.2(i)]{InterpolationCWHM}), the second inequality holds since, again given this normalisation, the embedding of $H^s(\Omega)=H^{\theta m}(\Omega)$ into $K_{\theta,2}((L^2(\Omega),H^m(\Omega)))$ has norm $\leq 1$ \cite[Lem.~4.2]{InterpolationCWHM},\footnote{Note that in general it does not hold that $H^{\theta m}(\Omega)=K_{\theta,2}((L^2(\Omega),H^m(\Omega)))$ - see, e.g., \cite[Lem.~4.8]{InterpolationCWHM}. }
and the final equality holds because $\theta m = s$. 

\subsection*{Proof of Corollary \ref{cor:Duality}}

We adopt a duality argument, similar to that used in the proofs of \cite[Thm.~4.1.33]{sauter-schwab11} and \cite[Lemma A.1]{BEMfract}. 
Let $-m\leq s_1\leq 0\leq s_2\leq m$ and let $u\in H^{s_2}(\Omega)\subset L^2(\Omega)$. By the duality between $\tH^{-s}(\Omega)$ and $\tH^s(\Omega)$ (see, e.g., \cite[Thm~3.3 \& \S3.1.3]{ChaHewMoi:13}), the $L^2$-orthogonality of the projection $\Pi$ (which we define as in the proof of Corollary \ref{cor:int} to be orthogonal projection onto $V_{h,p}$ in $L^2(\Omega)$), the Cauchy-Schwarz inequality, and two applications of \eqref{e:Est2}, we have that 
\begin{align*}
\label{}
\min_{v\in V_{h,p}}\|u-v\|_{\tH^{s_1}(\Omega)} \leq \|u-\Pi u\|_{\tH^{s_1}(\Omega)} 
&= \sup_{0\neq v\in H^{-s_1}(\Omega)}\frac{|(u-\Pi  u,v)_{L_2(\Omega)}|}{\|v\|_{H^{-s_1}(\Omega)}} \\
&= \sup_{0\neq v\in H^{-s_1}(\Omega)}\frac{|(u-\Pi  u,v-\Pi  v)_{L_2(\Omega)}|}{\|v\|_{H^{-s_1}(\Omega)}}\\
&\leq  \sup_{0\neq v\in H^{-s_1}(\Omega)}\frac{\|u-\Pi  u\|_{L_2(\Omega)}\|v-\Pi  v\|_{L_2(\Omega)}}{\|v\|_{H^{-s_1}(\Omega)}}\\
& \leq (C_{n,h_0,m})^{(s_2-s_1)/m}\left(\frac{h}{p+1}\right)^{s_2-s_1}\|u\|_{H^{s_2}(\Omega)}.
\end{align*}

\subsection*{Proof of Corollary \ref{cor:Duality2}}

Let $-m\leq s_1\leq 0$ and let $\Pi^{s_1}:\tH^{s_1}(\Omega)\to V_{h,p}$ denote orthogonal projection onto $\widetilde{V}_{h,p}$ in $\tH^{s_1}(\Omega)$. 
Applying \eqref{eq:Hsboundh_dn} with $s_2=0$ gives
\begin{align}
\label{eq:Hs1bound}
\|u-\Pi^{s_1} u\|_{\tH^{s_1}(\Omega)}   
\leq  
(C_{n,h_0,m})^{-s_1/m} \left(\frac{h}{p+1}\right)^{-s_1}
\|u\|_{\tH^{0}(\Omega)}, \qquad u\in\tH^{0}(\Omega),
\end{align}
and we
also have the trivial estimate
\begin{align}
\label{eq:Hs1boundtrivial}
\|u-\Pi^{s_1} u\|_{\tH^{s_1}(\Omega)}   
\leq  
\|u\|_{\tH^{s_1}(\Omega)}, \qquad u\in\tH^{s_1}(\Omega).
\end{align}
Then \eqref{eq:Hsboundh_dn2} follows for $-m\leq s_1\leq s_2\leq 0$ by applying interpolation to $I-\Pi^{s_1}$, 
using \eqref{eq:Hs1bound} and \eqref{eq:Hs1boundtrivial}.

\section*{Acknowledgements}
The author thanks Simon Chandler-Wilde, Emmanuil Georgoulis, Andrea Moiola and Iain Smears for helpful discussions in relation to this work. Part of this work was carried out while the author was in residence at Institut Mittag-Leffler in Djursholm, Sweden in autumn 2025, supported by the Swedish Research Council under grant no.~2021-06594. 

\bibliography{BEMbib_short2014} %

\begin{thebibliography}{10}

\bibitem{Achdou07}
{\sc Y.~Achdou, C.~Sabot, and N.~Tchou}, {\em Transparent boundary conditions
  for the {H}elmholtz equation in some ramified domains with a fractal
  boundary}, J. Comput. Phys., 220 (2007), pp.~712--739.

\bibitem{atkinson1997numerical}
{\sc K.~E. Atkinson}, {\em {The Numerical Solution of Integral Equations of the
  Second Kind}}, Cambridge University Press, 1997.

\bibitem{babuska1987optimal}
{\sc I.~Babu{\v{s}}ka and M.~Suri}, {\em The optimal convergence rate of the
  p-version of the finite element method}, SIAM J. Numer. Anal., 24 (1987),
  pp.~750--776.

\bibitem{babuska1987hp}
{\sc I.~Babu{\v{s}}ka and M.~Suri}, {\em {The $hp$ version of the finite
  element method with quasiuniform meshes}}, ESAIM: Mod{\'e}l. Math. Anal.
  Num{\'e}r., 21 (1987), pp.~199--238.

\bibitem{bagnerini2006finite}
{\sc P.~Bagnerini, A.~Buffa, and E.~Vacca}, {\em Finite elements for a
  prefractal transmission problem}, C. R. Math., 342 (2006), pp.~211--214.

\bibitem{Bagnerini13}
{\sc P.~Bagnerini, A.~Buffa, and E.~Vacca}, {\em Mesh generation and numerical
  analysis of a {G}alerkin method for highly conductive prefractal layers},
  Appl. Numer. Math., 65 (2013), pp.~63--78.

\bibitem{FractalTransmission}
{\sc J.~Bannister, D.~P. Hewett, and A.~Gibbs}, {\em {Acoustic scattering by
  fractal inhomogeneities via geometry-conforming Galerkin methods for the
  Lippmann-Schwinger equation}}.
\newblock arxiv.org/abs/2602.05005.

\bibitem{brenner2008mathematical}
{\sc S.~C. Brenner and L.~R. Scott}, {\em {The Mathematical Theory of Finite
  Element Methods}}, Springer, 2008.

\bibitem{HausdorffDomain}
{\sc A.~M. Caetano, S.~N. Chandler-Wilde, X.~Claeys, A.~Gibbs, D.~P. Hewett,
  and A.~Moiola}, {\em {Integral equation methods for acoustic scattering by
  fractals}}, Proc. Roy. Soc. A, 481 (2025), p.~20230650.

\bibitem{HausdorffBEM}
{\sc A.~M. Caetano, S.~N. Chandler-Wilde, A.~Gibbs, D.~P. Hewett, and
  A.~Moiola}, {\em {A Hausdorff-measure boundary element method for acoustic
  scattering by fractal screens}}, Numer. Math., 156 (2024), p.~463–532.

\bibitem{CaChHe25}
{\sc A.~M. Caetano, S.~N. Chandler-Wilde, and D.~P. Hewett}, {\em {Properties
  of IFS attractors with non-empty interiors and related rough sets, and
  associated function spaces and scattering problems}}.
\newblock arxiv.org/abs/2511.15213.

\bibitem{caetano2019density}
{\sc A.~M. Caetano, D.~P. Hewett, and A.~Moiola}, {\em {Density results for
  Sobolev, Besov and Triebel-Lizorkin spaces on rough sets}}, J. Funct. Anal.,
  281 (2021), p.~109019.

\bibitem{Calderon61}
{\sc A.~P. Calder{\'o}n}, {\em Lebesgue spaces of differentiable functions and
  distributions}, Proc. Symp. Pure Math., 4 (1961), pp.~33--49.

\bibitem{cangiani2022hp}
{\sc A.~Cangiani, Z.~Dong, and E.~Georgoulis}, {\em {$hp$-version discontinuous
  Galerkin methods on essentially arbitrarily-shaped elements}}, Math. Comp.,
  91 (2022), pp.~1--35.

\bibitem{capitanelli2015weighted}
{\sc R.~Capitanelli and M.~A. Vivaldi}, {\em Weighted estimates on fractal
  domains}, Mathematika, 61 (2015), pp.~370--384.

\bibitem{InterpolationCWHM}
{\sc S.~N. Chandler-Wilde, D.~P. Hewett, and A.~Moiola}, {\em Interpolation of
  {H}ilbert and {S}obolev spaces: quantitative estimates and counterexamples},
  Mathematika, 61 (2015), pp.~414--443.

\bibitem{ChaHewMoi:13}
\leavevmode\vrule height 2pt depth -1.6pt width 23pt, {\em {S}obolev spaces on
  non-{L}ipschitz subsets of $\mathbb{R}^n$ with application to boundary
  integral equations on fractal screens}, Integr. Equat. Operat. Th., 87
  (2017), pp.~179--224.

\bibitem{BEMfract}
{\sc S.~N. Chandler-Wilde, D.~P. Hewett, A.~Moiola, and J.~Besson}, {\em
  Boundary element methods for acoustic scattering by fractal screens}, Numer.
  Math., 147 (2021), pp.~785--837.

\bibitem{Ciarlet78}
{\sc P.~G. Ciarlet}, {\em The Finite Element Method for Elliptic Problems},
  North-Holland, 1978.

\bibitem{creo2022transmission}
{\sc S.~Creo, M.~R. Lancia, and P.~Vernole}, {\em {Transmission problems for
  the fractional p-Laplacian across fractal interfaces}}, Discrete Contin. Dyn.
  Syst. - Ser. S, 15 (2022).

\bibitem{NonDisjointQuad}
{\sc A.~Gibbs, D.~P. Hewett, and B.~Major}, {\em Numerical evaluation of
  singular integrals on non-disjoint self-similar fractal sets}, Numer. Alg.,
  97 (2024), pp.~311--343.

\bibitem{HausdorffQuadrature}
{\sc A.~Gibbs, D.~P. Hewett, and A.~Moiola}, {\em {Numerical evaluation of
  singular integrals on fractal sets}}, Numer. Alg., 92 (2023), pp.~2071--2124.

\bibitem{GoMoHe}
{\sc S.~G\'omez, D.~P. Hewett, and A.~Moiola}, {\em {A discontinuous Galerkin
  method with fractal elements}}.
\newblock arxiv.org/abs/2604.11093.

\bibitem{hinz2021non}
{\sc M.~Hinz, A.~Rozanova-Pierrat, and A.~Teplyaev}, {\em {Non-Lipschitz
  uniform domain shape optimization in linear acoustics}}, SIAM J. Control
  Optim., 59 (2021), pp.~1007--1032.

\bibitem{Jones}
{\sc P.~W. Jones}, {\em Quasiconformal mappings and extendability of functions
  in {S}obolev spaces}, Acta Math., 147 (1981), pp.~71--88.

\bibitem{JoWa84}
{\sc A.~Jonsson and H.~Wallin}, {\em Function {S}paces on {S}ubsets of
  {${\mathbb R}^n$}}, Math. Rep., 2 (1984).

\bibitem{Kress}
{\sc R.~Kress}, {\em Linear Integral Equations}, Springer, 2nd ed., 1999.

\bibitem{Rogers}
{\sc L.~G. Rogers}, {\em Degree-independent {S}obolev extension on locally
  uniform domains}, J. Funct. Anal., 235 (2006), pp.~619--665.

\bibitem{Rychkov}
{\sc V.~S. Rychkov}, {\em On restrictions and extensions of the {B}esov and
  {T}riebel--{L}izorkin spaces with respect to {L}ipschitz domains}, J. London
  Math. Soc., 60 (1999), pp.~237--257.

\bibitem{sauter-schwab11}
{\sc S.~A. Sauter and C.~Schwab}, {\em Boundary Element Methods}, Springer,
  2011.

\bibitem{Sc:98}
{\sc C.~Schwab}, {\em $p$- and $hp$-Finite Element Methods}, Oxford University
  Press, 1998.

\bibitem{smears2015}
{\sc I.~Smears}, {\em {Discontinuous Galerkin finite element approximation of
  Hamilton–Jacobi–Bellman equations with Cordes coefficients}}, PhD thesis,
  University of Oxford, 2015.

\bibitem{Stein}
{\sc E.~M. Stein}, {\em Singular Integrals and Differentiability Properties of
  Functions}, Princeton University Press, 1970.

\bibitem{Wolff1982}
{\sc T.~H. Wolff}, {\em A note on interpolation spaces}, in {Harmonic Analysis:
  Proceedings of a conference held at the University of Minnesota, Minneapolis,
  April 20--30, 1981}, Springer, 1982, pp.~199--204.

\end{thebibliography}
\bibliographystyle{siam}
\end{document}